\documentclass{amsart}

\usepackage{amsmath,amstext,amssymb,amsopn,amsthm,mathrsfs}

\usepackage{psfrag}
\usepackage[dvips]{color}
\usepackage{multicol}
\usepackage{picinpar}
\usepackage{enumerate}
\usepackage[hyphens]{url}        
\usepackage[breaklinks,colorlinks=true,linkcolor=blue,citecolor=blue, urlcolor=blue]{hyperref}                                                 
\usepackage[utf8x]{inputenc}
\usepackage{txfonts}
\usepackage{graphicx}
\usepackage[T1]{fontenc}
\usepackage{longtable}
\theoremstyle{plain}
\newtheorem{teo}{Theorem}[section]
\newtheorem{defi}{Definition}[section]
\newtheorem{corollary}{Corollary}[section]
\newtheorem{lemma}{Lemma}[section]

\newtheorem{obs}{Observation}[section]
\numberwithin{equation}{section}
\newcommand{\dem}{\medskip \par \noindent \mbox{\bf Proof. }}

\begin{document}

\title[The Boundedness of the Ornstein-Uhlenbeck semigroup on $L^{p(\cdot)}(\gamma_d).$ ] {The Boundedness of the Ornstein-Uhlenbeck semigroup on variable Lebesgue spaces with respect to the Gaussian measure.}

\author{Jorge Moreno}
\address{Departamento de Matem\'{a}tica,  Decanato de Ciencia y
Tecnologia, UCLA
 Apartado 400 Barquisimeto 3001 Venezuela.}
\email{jorge.moreno@ucla.edu.ve}
   \author{Ebner Pineda}
\address{Departamento de Matem\'{a}tica,  Facultad de Ciencias Naturales y Matem\'aticas, ESPOL Guayaquil 09-01-5863, Ecuador.}
\email{epineda@espol.edu.ec}
\author{Wilfredo~O.~Urbina}
\address{Department of Mathematics and Actuarial Sciences, Roosevelt University, Chicago, IL,
   60605, USA.}
\email{wurbinaromero@roosevelt.edu}



\begin{abstract}
The main result of this work is the proof of the boundedness of  the  Ornstein-Uhlenbeck  semigroup $ \{T_t \}_{t\geq 0} $ in $ {\mathbb R}^d $ on Gaussian variable Lebesgue spaces under a condition of regularity on $p(\cdot)$ following \cite{DalSco} and  \cite{TesSon}. As a consequence of this result, we obtain  the boundedness of Poisson-Hermite semigroup and the boundedness of the  Gaussian Bessel potentials of order $\beta> 0$.
\end{abstract}

\maketitle

\section{Introduction and Preliminaries}

The  Ornstein-Uhlenbeck semigroup $\{T_t\}_{t\geq 0} $ is the semigroup of operators generated in $L^2(\gamma_d)$ by the Ornstein-Uhlenbeck operator 
\begin{equation}
L =  \frac{1}{2}\Delta_x - \langle x, \nabla_x \rangle= \sum_{i=1}^d \Big[\frac{1}{2} \frac{\partial^2}{\partial x_i^2} - x_i \frac{\partial }{\partial x_i}\Big]
\end{equation}
 as infinitesimal generator, i.e. formally $T_t = e^{-tL}$. In view of the spectral theorem, for $f  = \sum_{k=0}^{\infty} {\bf J}_k f  \in L^2(\mathbb{R}^{d},\gamma_d)$ and $t \geq 0,$ $T_t$ is defined as 
\begin{equation}\label{repOUAbel}
T_t f =  \sum_{\nu} e^{-t|\nu|}  \langle f, {\vec h}_{\nu} \rangle_{\gamma_d} {\vec h}_{\nu}
=  \sum_{k=0}^{\infty}   e^{-tk} \sum_{|\nu| = k}  \langle f, 
{\vec h}_{\nu} \rangle_{\gamma_d} {\vec h}_{\nu} = \sum_{k=0}^{\infty} e^{-tk} {\bf J}_k f ,
\end{equation}
where  $\{\vec{h}_{\nu}\}_{\nu}$  are the normalized Hermite polynomials in $d$ variables, 
$${\bf J}_k f  =  \sum_{|\nu| = k}   \langle f, 
{\vec h}_{\nu} \rangle_{\gamma_d} {\vec h}_{\nu}$$
 is the orthogonal projection of $L^2 (\mathbb{R}^{d},\gamma_d)$ onto 
$$\mathcal{C}_k = \overline{\mbox{span}\left(\left\{ \vec{h}_{\nu}: |\nu | = k\right\}\right)}^{L^2 (\mathbb{R}^{d},\gamma_d)}.$$

It can be proved, using Mehler's formula, that the Ornstein-Uhlenbeck semigroup has an integral representation as
\begin{eqnarray}
\nonumber T_t f (x) &=& \frac{1}{(1-e^{-2t})^{d/2}} \int_{{\mathbb R}^d}  e^{- \frac{e^{-2t}(|y|^2+|x|^2) - 2 e^{-t} \langle x, y \rangle}{1-e^{-2t}}} f(y)\gamma_d(dy)\\
&=&  \frac{1}{\pi^{d/2}(1-e^{-2t})^{d/2}}\int_{{\mathbb R}^d} 
e^{- \frac{|y-e^{-t}x|^2}{1-e^{-2t}}} f(y) dy.
\end{eqnarray}
for all $f\in L^{1}(\mathbb{R}^{d},\gamma_{d})$.
Taking the change of variable $s=1-e^{-2t}$, we obtain that
\begin{equation}\label{OUAltrep}
  T_{t}f(x) =\frac{1}{(\pi s)^{d/2}}  \int_{{\mathbb R}^d} e^{-\frac{|y-\sqrt{(1-s)}x|^{2}}{s}}f(y)dy.
\end{equation}

It is well know that the  Ornstein-Uhlenbeck  semigroup $ \{T_t \}_{t\geq 0} $ in $ {\mathbb R}^d $ is a {\em Markov operator semigroup} in $ L^p (\mathbb{R}^{d},\gamma_d), 1 \leq p \leq \infty$,  i.e. a positive conservative symmetric diffusion semigroup, strongly $L^p$-continuous in $ L^p (\mathbb{R}^{d},\gamma_d), 1 \leq p \leq \infty$; with the Ornstein-Uhlenbeck operator $L$ as its infinitesimal generator, see \cite{ba2}, \cite{ba3} or \cite{urbina2019}. Its properties can be obtained directly from the general theory of Markov semigroups, see  \cite{ba3} or \cite{ur7}. \\

We prove that  the  Ornstein-Uhlenbeck  semigroup $ \{T_t \}_{t\geq 0} $ in $ {\mathbb R}^d $ is also bounded for  Gaussian variable Lebesgue spaces $L^{p(\cdot)}(\mathbb{R}^{d},\gamma_d),$ under a condition of regularity on $p(\cdot).$ 

\begin{teo}\label{acotacionsemigrupoOU}
Let $p(\cdot)\in\mathcal{P}_{\gamma_{d}}^{\infty}(\mathbb{R}^{d})\cap LH_{0}(\mathbb{R}^{d})$ with $1<p_{-}\leq p_{+}<\infty$. Then there exists a constant $C>0$ such that
\begin{eqnarray*}
  ||T_{t}f||_{p(\cdot),\gamma_{d}}\leq C||f||_{p(\cdot),\gamma_{d}}
\end{eqnarray*}
for all  $f\in L^{p(\cdot)}(\mathbb{R}^{d},\gamma_{d})$ and $t>0$.
\end{teo}

The maximal function of the Ornstein-Uhlenbeck semigroup is defined by
\begin{eqnarray*}
T^{*}f(x)=\sup_{t>0}T_{t}f(x),
\end{eqnarray*}
for all $x\in\mathbb{R}^{d}.$

As a consequence of the proofs of Theorems \ref{acotacionsemigrupoOUlocal} and \ref{acotacionsemigrupoOUglobal} we obtain,
\begin{corollary}\label{operadormaxOU}
Let $p(\cdot)\in\mathcal{P}_{\gamma_{d}}^{\infty}(\mathbb{R}^{d})\cap LH_{0}(\mathbb{R}^{d})$ with $1<p_{-}\leq p_{+}<\infty.$ Then there exists a constant $C>0$ such that
\begin{eqnarray*}
  ||T^{*}f||_{p(\cdot),\gamma_{d}}\leq C||f||_{p(\cdot),\gamma_{d}}
\end{eqnarray*}
for all  $f\in L^{p(\cdot)}(\mathbb{R}^{d},\gamma_{d})$.\\
\end{corollary}

Two important remarks are needed here. First, observe that in Theorem \ref{acotacionsemigrupoOU} we can not conclude, as in the classical case that $T_t$ is a contraction. We do not know if that is actually true for this case, or if simply a problem of the method of proof. Therefore questions like some form of hypercontractivity for the semigroup  in this context are totally unknown. Second, the method of the proof give us immediately the boundedness of the maximal function of the Ornstein-Uhlenbeck semigroup.\\

Additionally, let us consider the {\em Poisson-Hermite semigroup} as the semigroup subordinated to the Ornstein-Uhlenbeck semigroup,  using the {\em Bochner's subordination formula}, see E. Stein \cite{ste70},
defined then as,
\begin{eqnarray}\label{Poisson-HermRep1}
P_t f(x) & = & \frac{1}{\sqrt \pi} \int_0^{\infty} \frac{e^{-u}}{\sqrt u} T_{(t^2/4u)}f(x) du
\nonumber \\
	& = &  \frac{1}{2\pi^{(d+1)/2}} \int_{{\mathbb R}^d} \int_0^1 t \frac{\exp(t
^2/4\log r)}{(-\log r)^{3/2}} 
 \frac{\exp\Big(\frac{-|y-rx|^2}{1-r^2}\Big)}{(1-r^2)^{d/2}} \frac{dr}{r} f(y) dy. 
\end{eqnarray}

It is also well known, that the Poisson-Hermite semigroup   $\{ P_t \}_{t \geq 0}$ is a  strongly continuous,  symmetric, conservative semigroup of positive contractions in 
 $L^p(\gamma_d)$, $ 1 \leq p < \infty$, with infinitesimal generator $(-L)^{1/2}$. As a consequence of the boundedness of $\{T_{t}\},$  we will prove that  it is also bounded for  Gaussian variable Lebesgue spaces $L^{p(\cdot)}(\mathbb{R}^{d},\gamma_d)$ under the same condition of regularity on $p(\cdot).$ \\
 \begin{teo}\label{acotacionsemigrupoPH}
Let $p(\cdot)\in\mathcal{P}_{\gamma_{d}}^{\infty}(\mathbb{R}^{d})\cap LH_{0}(\mathbb{R}^{d})$ with $1<p_{-}\leq p_{+}<\infty.$ Then there exists a constant $C>0$ such that
\begin{eqnarray*}
  ||P_{t}f||_{p(\cdot),\gamma_{d}}\leq C||f||_{p(\cdot),\gamma_{d}}
\end{eqnarray*}
for all  $f\in L^{p(\cdot)}(\mathbb{R}^{d},\gamma_{d})$ and $t>0$.
\end{teo}

 Finally, the Gaussian Bessel potential of order $\beta>0$, $\mathscr{J}_{\beta}$ is defined as
\begin{equation}
\mathscr{J}_{\beta}f(x)=\frac{1}{\Gamma(\beta)}\int_{0}^{+\infty}s^{\beta-1}e^{-s}P_{s}f(x)\; ds
\end{equation}
for all $x\in\mathbb{R}^{d}.$

It can be proved, using P. A. Meyer's multiplier theorem, that the Gaussian Bessel potentials $\mathscr{J}_\beta$ are $L^p(\gamma_d)$-bounded $1 < p < \infty.$  Moreover we will see that as consequence of Theorem \ref{acotacionsemigrupoPH} we obtain the boundedness of Gaussian Bessel potential on $L^{p(\cdot)}(\mathbb{R}^{d},\gamma_d).$

  \begin{teo}\label{acotacionpotencialBessel}
Let $p(\cdot)\in\mathcal{P}_{\gamma_{d}}^{\infty}(\mathbb{R}^{d})\cap LH_{0}(\mathbb{R}^{d})$ with $1<p_{-}\leq p_{+}<\infty.$ Then there exists a constant $C>0$ such that
\begin{eqnarray*}
  ||\mathscr{J}_{\beta}f||_{p(\cdot),\gamma_{d}}\leq C||f||_{p(\cdot),\gamma_{d}}
\end{eqnarray*}
for all $f\in L^{p(\cdot)}(\mathbb{R}^{d},\gamma_{d})$ and $\beta>0$.\\
\end{teo}

Now, for completeness, let us get more background on variable Lebesgue spaces with respect to a Borel measure $\mu$.

Any $\mu$-measurable function $p(\cdot):\mathbb{R}^{d}\rightarrow [1,\infty]$ is an exponent function; the set of all the exponent functions will be denoted by  $\mathcal{P}(\mathbb{R}^{d},\mu)$. For $E\subset\mathbb{R}^{d}$ we set $$p_{-}(E)=\text{ess}\inf_{x\in E}p(x) \;\text{and}\; p_{+}(E)=\text{ess}\sup_{x\in E}p(x).$$
We use the abbreviations $p_{+}=p_{+}(\mathbb{R}^{d})$ and $p_{-}=p_{-}(\mathbb{R}^{d})$.

\begin{defi}\label{deflogholder}
Let $E\subset \mathbb{R}^{d}$. We say that $\alpha(\cdot):E\rightarrow\mathbb{R}$ is locally log-H\"{o}lder continuous, and denote this by $\alpha(\cdot)\in LH_{0}(\mathbb{R}^{d})$, if there exists a constant $C_{1}>0$ such that
			\begin{eqnarray*}
				|\alpha(x)-\alpha(y)|&\leq&\frac{C_{1}}{log(e+\frac{1}{|x-y|})}
			\end{eqnarray*}
			for all $x,y\in E$. We say that $\alpha(\cdot)$ is log-H\"{o}lder continuous at infinity with base point at $x_{0}\in \mathbb{R}^{d}$, and denote this by $\alpha(\cdot)\in LH_{\infty}(\mathbb{R}^{d})$, if there exist  constants $\alpha_{\infty}\in\mathbb{R}$ and $C_{2}>0$ such that
			\begin{eqnarray*}
				|\alpha(x)-\alpha_{\infty}|&\leq&\frac{C_{2}}{log(e+|x-x_{0}|)}
			\end{eqnarray*}
			for all $x\in E$. We say that $\alpha(\cdot)$ is log-H\"{o}lder continuous, and denote this by $\alpha(\cdot)\in LH(\mathbb{R}^{d})$ if both conditions are satisfied.
			The maximum, $\max\{C_{1},C_{2}\}$ is called the log-H\"{o}lder constant of $\alpha(\cdot)$.
\end{defi}

\begin{defi}\label{defPdlog}
			We denote $p(\cdot)\in\mathcal{P}_{d}^{log}(\mathbb{R}^{d})$, if $\frac{1}{p(\cdot)}$ is log-H\"{o}lder continuous and  denote by $C_{log}(p)$ or $C_{log}$ the log-H\"{o}lder constant of $\frac{1}{p(\cdot)}$.
		\end{defi}

We will need the following technical result, for its proof see Lemma 3.26 in \cite{dcruz}.			
\begin{lemma}\label{lema3.26CU}
Let $\rho(\cdot):\mathbb{R}^{d}\rightarrow[0,\infty)$ be such that $\rho(\cdot)\in LH_{\infty}(\mathbb{R}^{d})$, $0<\rho_{\infty}<\infty$, and let $R(x)=(e+|x|)^{-N}$, $N>d/\rho_{-}$. Then there exists a constant $C$ depending on $d$, $N$ and the $LH_{\infty}$ constant of $r(\cdot)$ such that given any set $E$ and
any function $F$ with $0\leq F(y)\leq 1$ for $y\in E$,
\begin{eqnarray}
  \int_{E}F^{\rho(y)}(y)dy &\leq& C\int_{E}F(y)^{\rho_{\infty}}dy + \int_{E}R^{\rho_{-}}(y)dy,\label{3.26.1} \\
  \int_{E}F^{\rho_{\infty}}(y)dy &\leq& C\int_{E}F^{r(y)}(y)dy + \int_{E}R^{\rho^{-}}(y)dy.\label{3.26.2}
\end{eqnarray}
\end{lemma}

\begin{defi}
For a $\mu$-measurable function $f:\mathbb{R}^{d}\rightarrow \mathbb{R}$, we define the modular \begin{equation}
\rho_{p(\cdot),\mu}(f)=\displaystyle\int_{\mathbb{R}^{d}\setminus\Omega_{\infty}}|f(x)|^{p(x)}\mu(dx)+\|f\|_{L^{\infty}(\Omega_{\infty},\mu)},
\end{equation}
and the norm
\begin{equation}
\|f\|_{L^{p(\cdot)}(\mathbb{R}^{d},\mu)}=\inf\left\{\lambda>0:\rho_{p(\cdot),\mu}(f/\lambda)\leq 1\right\}.
\end{equation}

\end{defi}

\begin{defi}The variable exponent Lebesgue space on $\mathbb{R}^{d}$, $L^{p(\cdot)}(\mathbb{R}^{d},\mu)$ consists on those $\mu\_$measurable functions $f$ for which there exists $\lambda>0$ such that $\rho_{p(\cdot),\mu}\left(\frac{f}{\lambda}\right)<\infty,$ i.e.
\begin{equation*}
L^{p(\cdot)}(\mathbb{R}^{d},\mu) =\left\{f:\mathbb{R}^{d}\to \mathbb{R}: f \; \text{measurable } \; \rho_{p(\cdot),\mu}\left(\frac{f}{\lambda}\right)<\infty, \; \text{for some} \;\lambda>0\right\}.
\end{equation*}

\end{defi}

If $\mathcal{B}$ is a family of balls (or cubes) in $\mathbb{R}^{d}$, we say that $\mathcal {B}$ is $N$-finite if it has bounded overlappings for $N$, this is $\displaystyle\sum_{B\in\mathcal{B}}\chi_{B}(x)\leq N$ for all $x\in\mathbb {R}^{d}$; in other words, there is only $N$ balls (resp cubes) that intersect at the same time.\\

The following definition was introduced for the first time by Berezhno\v{\i} in \cite{Berez}, defined for family of disjoint balls or cubes. In the context of variable spaces, it has been considered in \cite{LibroDenHarjHas}, allowing the family to have bounded overlappings.\\
\begin{defi}
Given an exponent $p(\cdot)\in\mathcal{P}(\mathbb{R}^{d})$, we will say that $p(\cdot)\in\mathcal{G}$, if for every family of balls (or cubes) $\mathcal{B}$ which is $N$-finite,
\begin{eqnarray*}
  \sum_{B\in\mathcal{B}}||f\chi_{B}||_{p(\cdot)}||g\chi_{B}||_{p^{'}(\cdot)} &\lesssim& ||f||_{p(\cdot)}||g||_{p^{'}(\cdot)}
\end{eqnarray*}
for all functions $f\in L^{p(\cdot)}(\mathbb{R}^{d})$ and $g\in L^{p^{'}(\cdot)}(\mathbb{R}^{d})$. The constant only depends on N.
\end{defi}

\begin{lemma}[Teorema 7.3.22 in \cite{LibroDenHarjHas}]\label{implication1}
If $p(\cdot)\in LH(\mathbb{R}^{d})$, then $p(\cdot)\in\mathcal{G}$
\end{lemma}

As usual in what follows $C$ represents a constant that is not necessarily the same in each occurrence; also we will used the notation: given two functions $f$, $g$, the symbols $\lesssim$ and $\gtrsim$ denote, that there is a constant $c$ such that $f\leq cg$ and $cf\geq g$, respectively. When both inequalities are satisfied, that is, $f\lesssim g\lesssim f$, we will denote $f\approx g$. 

\section{Proofs of the main results.}

In this section we are going to consider Lebesgue variable spaces with respect to the Gaussian measure $\gamma_d,$ $L^{p(\cdot)}(\mathbb{R}^{d},\gamma_d).$  The next condition was introduced by E. Dalmasso and R. Scotto in \cite{DalSco}.

\begin{defi}\label{defipgamma}
Let $p(\cdot)\in\mathcal{P}(\mathbb{R}^{d},\gamma_{d})$, we say that $p(\cdot)\in\mathcal{P}_{\gamma_{d}}^{\infty}(\mathbb{R}^{d})$ if there exist constants $C_{\gamma_{d}}>0$ and $p_{\infty}\geq1$ such that
\begin{equation}
   |p(x)-p_{\infty}|\leq\frac{C_{\gamma_{d}}}{|x|^{2}},
\end{equation}
for $x\in\mathbb{R}^{d}\setminus\{(0,0,\ldots,0)\}.$
\end{defi}

\begin{obs}\label{obs4.1}
If $p(\cdot)\in\mathcal{P}_{\gamma_{d}}^{\infty}(\mathbb{R}^{d})$, then $p(\cdot)\in LH_{\infty}(\mathbb{R}^{d})$
\end{obs}
\begin{lemma}\label{lemaequiPgamma}
If $1<p_{-}\leq p_{+}<\infty,$ the following statements are equivalent
\begin{itemize}
  \item [(i)] $p(\cdot)\in\mathcal{P}_{\gamma_{d}}^{\infty}(\mathbb{R}^{d})$
  \item [(ii)] There exists $p_{\infty}>1$ such that
  \begin{eqnarray}
    C_{1}^{-1}\leq e^{-|x|^{2}(p(x)/p_{\infty}-1)}\leq C_{1} &\;\;\hbox{and}\;\;& C_{2}^{-1}\leq e^{-|x|^{2}(p^{'}(x)/p^{'}_{\infty}-1)}\leq C_{2},
  \end{eqnarray}
  for all $x\in\mathbb{R}^{d}$, where $C_{1}=e^{C_{\gamma_{d}}/p_{\infty}}$ and $C_{2}=e^{C_{\gamma_{d}}(p_{-})^{'}/p_{\infty}}$.
\end{itemize}
\end{lemma}

 Definition \ref{defipgamma} with observation \ref{obs4.1} and Lemma \ref{lemaequiPgamma} end up strengthening the regularity conditions on the exponent functions $p(\cdot)$ to obtain the boundedness of the semigroup $\{T_{t}\}$. As a consequence of Lemma \ref{implication1}, we have
\begin{corollary}\label{solapamientoacotadoG}
If $p(\cdot)\in\mathcal{P}_{\gamma_{d}}^{\infty}(\mathbb{R}^{d})\cap LH_{0}(\mathbb{R}^{d})$, then $p(\cdot)\in\mathcal{G}$
\end{corollary}

\begin{lemma}\label{equivnormp}
Let $p(\cdot)\in\mathcal{P}(\mathbb{R}^{d},\gamma_{d})$, then
\begin{eqnarray*}
\|f\|_{p(\cdot),\gamma_{d}}&\approx&\|fe^{-|\cdot|^{2}/p(\cdot)}\|_{p(\cdot)}\
\end{eqnarray*}
\end{lemma}
\begin{dem}
Let $A=\left\{\lambda>0:\;\rho_{p(\cdot)}\left(\frac{fe^{-|\cdot|^{2}/p(\cdot)}}{\lambda}\right)\leq1\right\}$ and $B=\left\{\lambda>0:\;\rho_{p(\cdot),\gamma_{d}}\left(\frac{f}{\lambda}\right)\leq1\right\}$. We will prove that $\inf(A)\lesssim\inf(B)$ and $\inf(B)\lesssim\inf(A)$. In fact, taking $\lambda\in A$ then 
$$\rho_{p(\cdot)}\left(\frac{fe^{-|\cdot|^{2}/p(\cdot)}}{\lambda}\right)=\int_{\mathbb{R}^{d}}\left| \frac{f(x)}{\lambda}\right|^{p(x)}e^{-|x|^{2}}dx\leq1$$
which implies
$$\rho_{p(\cdot),\gamma_{d}}\left(\frac{f}{\lambda}\right)\leq\frac{1}{\pi^{d/2}}\leq1$$
and then $\lambda\in B.$
Therefore, $A \subset B,$ and then
$\inf{B}\leq\inf{A}.$

On the other hand, taking $\lambda\in B$ then
$$\rho_{p(\cdot),\gamma_{d}}\left(\frac{f}{\lambda}\right)=\int_{\mathbb{R}^{d}}\left| \frac{f(x)}{\lambda}\right|^{p(x)}\frac{e^{-|x|^{2}}}{\pi^{d/2}}dx\leq1$$
which implies
$$\int_{\mathbb{R}^{d}}\left| \frac{f(x)e^{-|x|^{2}/p(x)}}{\lambda\pi^{d/2}}\right|^{p(x)}dx=\rho_{p(\cdot)}\left(\frac{fe^{-|\cdot|^{2}/p(\cdot)}}{\lambda\pi^{d/2}}\right)\leq1$$
and therefore $\lambda\in \pi^{-d/2}A$. Thus $\inf{A}\leq\pi^{d/2}\inf{B},$ and then $\inf(A)\approx\inf(B)$ 

Hence, we get 
\begin{eqnarray*}
\|f\|_{p(\cdot),\gamma_{d}}\approx\|fe^{-|\cdot|^{2}/p(\cdot)}\|_{p(\cdot)}. \qed
\end{eqnarray*}
\end{dem}

\subsection{Boundedness of the Ornstein-Uhlenbeck semigroup with the condition $\mathcal{P}^{\infty}_{\gamma_{d}}(\mathbb{R}^{d})$}

For $x \in \mathbb{R}^d$ let us consider admissible (or hyperbolic) balls, 
\begin{equation}
 B_{h}(x)= \{y\in\mathbb{R}^{d}:|x-y|\leq d(1\wedge1/|x|)\}.
 \end{equation}
 It is well known that the Gaussian measure is essentially constant on $ B_{h}(x),$ see \cite[Chapter 1]{urbina2019}.
 
 As it is nowadays standard in Gaussian harmonic analysis, we split $T_t$ into  {\em local part} and {\em global part}, using the change of variable $s=1-e^{-2t},$
  $$  T_{t}f(x)= T^{0}_{s}f(x)+T^{1}_{s}f(x),$$
  for $x\in\mathbb{R}^{d},$
  where
 $$
  T^{0}_{s}f(x):=\int_{B_{h}(x)}\frac{e^{-\frac{|\sqrt{(1-s)}x-y|^{2}}{s}}}{(\pi s)^{d/2}}f(y)dy, $$
  the local part, whichis the restriction of $T_t$  to the admissible ball $B_{h}(x)$ and
   $$T^{1}_{s}f(x):=\int_{B_{h}^{c}(x)}\frac{e^{-\frac{|\sqrt{(1-s)}x-y|^{2}}{s}}}{(\pi s)^{d/2}}f(y)dy$$
 the global part,  which is the restriction of $T_t$  to the complement of admissible ball $B_{h}(x).$

 Next, we will need the following technical lemma to handle the proof of boundedness of the local part, for the proof see \cite{urbina2019}, for an earlier formulation see also \cite{FabGutSco}.

\begin{lemma}\label{lemacubribolasB}
Let us define the secuence $x_{k}=\sqrt{k}$ for $k\in\mathbb{N}$. For this strictly increasing secuence, we obtain a family of disjoint balls $B_{j}^{k}$, for $k\in\mathbb{N}$ and $1\leq j\leq N_{k}$ with the following properties
\begin{description}
  \item[(i)] If $\tilde{B_{j}^{k}}=2B_{j}^{k}$, the colection $\mathcal{F}=\{B(0,1),\{\tilde{B_{j}^{k}}\}_{j,k}\}$ is a covering of $\mathbb{R}^{d}$;
  \item[(ii)] $\mathcal{F}$ has bounded overlappings;
  \item[(iii)] The center $y_{j}^{k}$ of $B_{j}^{k}$, satisfies $|y_{j}^{k}|=(x_{k+1}+x_{k})/2$;
  \item[(iv)] $diam(B_{j}^{k})=x_{k+1}-x_{k}=1/(2|y_{j}^{k}|);$
  \item[(v)] For all ball $B\in\mathcal{F}$, and all $x,y\in B$, $\gamma_{d}(x)\approx\gamma_{d}(y)$ with constants independents on $B$;
  \item[(vi)] There exists a uniform constant, $C_{n}>0$, such that, if $x\in B\in\mathcal{F}$, then $B_{h}(x)\subset C_{n}B:=\hat{B}$. Moreover, the colection $\hat{\mathcal{F}}=\{\hat{B}\}_{B\in\tilde{F}}$, also verifies properties (ii)-(v).
  \end{description}
\end{lemma}
Now, we present the boundedness of the local part of the semigroup $\{T_{t}\}$.
\begin{teo}\label{acotacionsemigrupoOUlocal}
Let $p(\cdot)\in\mathcal{P}_{\gamma_{d}}^{\infty}(\mathbb{R}^{d})\cap LH_{0}(\mathbb{R}^{d})$ with $1<p_{-}\leq p_{+}<\infty,$. There exists a constant $C>0$ such that
\begin{eqnarray*}
  ||T^{0}_{s}f||_{p(\cdot),\gamma_{d}}\leq C||f||_{p(\cdot),\gamma_{d}}
\end{eqnarray*}
for all function $f\in L^{p(\cdot)}(\mathbb{R}^{d},\gamma_{d}).$
\end{teo}

\begin{dem}
We follow the proof of Theorem 3.3. in \cite{DalSco}. Without lost of generality suppose that $f\geq0$.
\begin{eqnarray*}
T^{0}_{s}f(x)&:=&\int_{B_{h}(x)}\frac{e^{-\frac{|\sqrt{(1-s)}x-y|^{2}}{s}}}{(\pi s)^{d/2}}f(y)dy=\int_{\mathbb{R}^{d}}\frac{e^{-u(s)}}{(\pi s)^{d/2}}f(y)\chi_{B_{h}(x)}(y)dy\\
&=&\int_{\mathbb{R}^{d}}M(s,x,y)f(y)\chi_{B_{h}(x)}(y)dy
\end{eqnarray*}
where $M(s,x,y)=\displaystyle\frac{e^{-u(s)}}{(\pi s)^{d/2}}$,   $u(s)=\displaystyle\frac{|\sqrt{(1-s)}x-y|^{2}}{s}$ and\\
 $B_{h}(x)= \{y\in\mathbb{R}^{d}:|x-y|\leq d(1\wedge1/|x|)\}.$\\\\
Following \cite{TesSon} we obtain that if $y\in B_{h}(x)$ then $e^{-u(s)}\leq C_{d}e^{-\frac{|x-y|^{2}}{s}}$
and therefore $M(s,x,y)\leq C_{d}\frac{e^{-\frac{|x-y|^{2}}{s}}}{(\pi s)^{d/2}}.$

Now, given $x\in \mathbb{R}^{d}$, by Lemma \ref{lemacubribolasB}, there exists $B\in\mathcal{F}$ such that $x\in B$ and $B_{h}(x)\subset\hat{B}$, in consequence,
\begin{eqnarray*}
\int_{\mathbb{R}^{d}}M(s,x,y)f(y)\chi_{B_{h}(x)}(y)dy&\leq&\int_{B_{h}(x)}C_{d}\frac{e^{-\frac{|x-y|^{2}}{s}}}{(\pi s)^{d/2}}f(y)dy\\
&\leq&C_{d}\int_{\mathbb{R}^{d}}\frac{e^{-\frac{|x-y|^{2}}{s}}}{s^{d/2}}f(y)\chi_{\hat{B}}(y)dy
\end{eqnarray*}

Set $\phi_{s}(z)=\frac{e^{-\frac{|z|^{2}}{s}}}{s^{d/2}},$ since $\{\phi_{s}\}_{s>0}$ is an approximation of  identity, we have

\begin{eqnarray*}
\int_{\mathbb{R}^{d}}M(s,x,y)f(y)\chi_{B_{h}(x)}(y)dy&\leq& C_{d}|(\phi_{s}*f\chi_{\hat{B}})(x)|\leq C_{d}M_{H-L}(f\chi_{\hat{B}})(x),
\end{eqnarray*}
if $x\in B.$ Hence
\begin{eqnarray*}
T^{0}_{s}f(x)\lesssim M_{H-L}(f\chi_{\hat{B}})(x)&=&M_{H-L}(f\chi_{\hat{B}})(x)\chi_{B}(x),
\end{eqnarray*}
 and therefore,
\begin{equation}\label{sumaximal}
T^{0}_{s}f(x)\lesssim\sum_{B\in\mathcal{F}}M_{H-L}(f\chi_{\hat{B}})(x)\chi_{B}(x),
\end{equation}
 for all $x\in \mathbb{R}^{d}.$
 
Let $f\in L^{p(\cdot)}(\mathbb{R}^{d},\gamma_{d})$. Using the characterization of norm by duality, 
\begin{eqnarray*}
  ||T^{0}_{s}f||_{p(\cdot),\gamma_{d}} &\leq& 2\sup_{||g||_{p'(\cdot),\gamma_{d}}\leq1}\int_{\mathbb{R}^{d}}T^{0}_{s}f(x)|g(x)|\gamma_{d}(dx)
\end{eqnarray*}
from (\ref{sumaximal}) and following again \cite{DalSco} we obtain that
\begin{eqnarray*}
  \int_{\mathbb{R}^{d}}M^{0}_{s}f(x)|g(x)|\gamma_{d}(dx) 
  &\lesssim& \sum_{B\in\mathcal{F}}e^{-|c_{B}|^{2}}\int_{\mathbb{R}^{d}}M_{H-L}(f\chi_{\hat{B}})(x)|g(x)|\chi_{B}(x)dx
\end{eqnarray*}
where $c_{B}$ is the center of the balls $B$ and $\hat{B}$.\\
Applying the H\"{o}lder's inequality and the boundedness of the maximal operator $M_{H-L}$ on $L^{p(\cdot)}(\mathbb{R}^{d})$, we get
\begin{eqnarray*}
 \int_{\mathbb{R}^{d}}T^{0}_{s}f(x)|g(x)|\gamma_{d}(dx)
 &\lesssim& \sum_{B\in\mathcal{F}}e^{-|c_{B}|^{2}\left(\frac{1}{p_{\infty}}+\frac{1}{p'_{\infty}}\right)}\|M_{H-L}(f\chi_{\hat{B}})\|_{p(\cdot)}\|g\chi_{B}\|_{p'(\cdot)}\\
   &\lesssim& \sum_{B\in\mathcal{F}}e^{-\frac{|c_{B}|^{2}}{p_{\infty}}}\|f\chi_{\hat{B}}\|_{p(\cdot)}e^{-\frac{|c_{B}|^{2}}{p'_{\infty}}}\|g\|_{p'(\cdot)}
\end{eqnarray*}
since $p(\cdot)\in\mathcal{P}_{\gamma_{d}}^{\infty}(\mathbb{R}^{d})$, by Lemma \ref{lemaequiPgamma}, we obtain that

%
%
%
%

\begin{eqnarray*}
e^{-\frac{|c_{B}|^{2}}{p_{\infty}}}\|f\chi_{\hat{B}}\|_{p(\cdot)}&\lesssim&\|f\chi_{\hat{B}}\|_{p(\cdot),\gamma_{d}}
\end{eqnarray*}
and
\begin{eqnarray*}
e^{-\frac{|c_{B}|^{2}}{p'_{\infty}}}\|g\chi_{\hat{B}}\|_{p'(\cdot)}&\lesssim&\|g\chi_{\hat{B}}\|_{p'(\cdot),\gamma_{d}}
\end{eqnarray*}
and by Lemma \ref{equivnormp}, we have that
\begin{eqnarray*}
\|f\chi_{\hat{B}}\|_{p(\cdot),\gamma_{d}}\approx\|f\chi_{\hat{B}}e^{-|\cdot|^{2}/p(\cdot)}\|_{p(\cdot)}\;&\hbox{and}&\;\|g\chi_{\hat{B}}\|_{p'(\cdot),\gamma_{d}}\approx\|g\chi_{\hat{B}}e^{-|\cdot|^{2}/p'(\cdot)}\|_{p'(\cdot)}
\end{eqnarray*}
therefore,

 \begin{eqnarray*}
 \int_{\mathbb{R}^{d}}T^{0}_{s}f(x)|g(x)|\gamma_{d}(dx)&\lesssim& \sum_{B\in\mathcal{F}}\|f\chi_{\hat{B}}e^{-|\cdot|^{2}/p(\cdot)}\|_{p(\cdot)}\|g\chi_{\hat{B}}e^{-|\cdot|^{2}/p'(\cdot)}\|_{p'(\cdot)}.
 \end{eqnarray*}
  Since the family of balls $\hat{\mathcal{F}}$ has bounded overlaps; applying Corollary \ref{solapamientoacotadoG}, to the functions $fe^{-|\cdot|^{2}/p(\cdot)}\in L^{p(\cdot)}(\mathbb{R}^{d})$ and $ge^{-|\cdot|^{2}/p'(\cdot)}\in L^{p'(\cdot)}(\mathbb{R}^{d})$ and again applying Lemma \ref{equivnormp}, we get

\begin{eqnarray*}
 \int_{\mathbb{R}^{d}}T^{0}_{s}f(x)|g(x)|\gamma_{d}(dx)&\lesssim&\|f\|_{p(\cdot),\gamma_{d}}\|g\|_{p'(\cdot),\gamma_{d}}.
 \end{eqnarray*}

 Taking supremum on all the functions $g\in L^{p'(\cdot)}(\mathbb{R}^{d},\gamma_{d})$ with\\
  $\|g\|_{p'(\cdot),\gamma_{d}}\leq1$, we obtain that
 \begin{eqnarray*}
 ||T^{0}_{s}f||_{p(\cdot),\gamma_{d}} &\lesssim& \sup_{||g||_{p'(\cdot),\gamma_{d}}\leq1}\int_{\mathbb{R}^{d}}M^{0}_{s}f(x)|g(x)|\gamma_{d}(dx)\\
 &\lesssim& \sup_{||g||_{p'(\cdot),\gamma_{d}}\leq1}\|f\|_{p(\cdot),\gamma_{d}}\|g\|_{p'(\cdot),\gamma_{d}}= \|f\|_{p(\cdot),\gamma_{d}}.
 \end{eqnarray*}
$\hfill\Box$

Finally, we will obtain the boundedness of the global part.
\begin{teo}\label{acotacionsemigrupoOUglobal}
Let $p(\cdot)\in\mathcal{P}_{\gamma_{d}}^{\infty}(\mathbb{R}^{d})\cap LH_{0}(\mathbb{R}^{d})$ con $1<p_{-}\leq p_{+}<\infty$. Then there exists a constant $C>0$ such that
\begin{eqnarray*}
  ||T^{1}_{s}f||_{p(\cdot),\gamma_{d}}\leq C||f||_{p(\cdot),\gamma_{d}}
\end{eqnarray*}
for all the function $f\in L^{p(\cdot)}(\mathbb{R}^{d},\gamma_{d}).$
\end{teo}
\begin{flushleft}
\dem
\end{flushleft}
Suppose that $f\geq0$. Again, we follow the proof of Theorem 3.5 in \cite{DalSco}.
\begin{eqnarray*}
  T^{1}_{s}f(x):=\int_{B^{c}(x)}\frac{e^{-\frac{|\sqrt{(1-s)}x-y|^{2}}{s}}}{(\pi s)^{d/2}}f(y)dy&=&\int_{B^{c}(x)}M(s,x,y)f(y)dy
\end{eqnarray*}
For  $x\in\mathbb{R}^{d}$ fix, set $E_{x}=\{y: b(x,y)>0\}$ where $b:=b(x,y)=2\left\langle x,y\right\rangle$. Given $y\in B^{c}(x)$, the following inequalities are satisfied:
\begin{enumerate}
  \item[(i)] If $b\leq0$, then
  \begin{equation}\label{acotacionMb-}
    M(s,x,y)\lesssim e^{-|y|^{2}}
  \end{equation}
  \item[(ii)] If $b>0$, then
  \begin{equation}\label{acotacionMb+}
    M(s,x,y)\lesssim \frac{e^{-u_{0}}}{t_{0}^{d/2}}
  \end{equation}
\end{enumerate}
where $a=|x|^{2}+|y|^{2}$, $t_{0}=2\sqrt{a^{2}-b^{2}}/(a+\sqrt{a^{2}-b^{2}})$ and $u_{0}=\frac{1}{2}(|y|^{2}-|x|^{2}+|x+y| |x-y|)$. For details see  \cite{TesSon} or \cite[Chapter 4]{urbina2019}.\\

Let $f\in L^{p(\cdot)}(\mathbb{R}^{d},\gamma_{d})$ with $\|f\|_{p(\cdot),\gamma_{d}}=1$. If $b\leq0$, applying (\ref{acotacionMb-}) and the H\"{o}lder's inequality for the exponent $p_{-}$ we obtain that
\begin{eqnarray*}
  I &=& \int_{\mathbb{R}^{d}}\left(\int_{B^{c}(x)\cap E^{c}_{x}}M(s,x,y)f(y)dy\right)^{p(x)}\gamma_{d}(dx) \\
   &\lesssim&  \int_{\mathbb{R}^{d}}\left(\int_{\mathbb{R}^{d}}(f(y))^{p_{-}}e^{-|y|^{2}}dy\right)^{p(x)/p_{-}}\gamma_{d}(dx).
\end{eqnarray*}
Moreover, $\rho_{p(\cdot),\gamma_{d}}(f)\leq 1$, implies that,
\begin{eqnarray*}
  I &\lesssim& 
   \int_{\mathbb{R}^{d}}\left(\int_{|f|>1}(f(y))^{p_{-}}e^{-|y|^{2}}dy+\int_{|f|\leq1}(f(y))^{p_{-}}e^{-|y|^{2}}dy\right)^{p(x)/p_{-}}\gamma_{d}(dx)\\
  &\lesssim& \int_{\mathbb{R}^{d}}\left(\int_{\mathbb{R}^{d}}(f(y))^{p(y)}\gamma_{d}(dy)+\int_{\mathbb{R}^{d}}\gamma_{d}(dy)\right)^{p(x)/p_{-}}\gamma_{d}(dx)\\
    &\lesssim&\int_{\mathbb{R}^{d}}(2)^{p(x)/p_{-}}\gamma_{d}(dx)= C_{d,p}.
\end{eqnarray*}
With this we obtain that $\|T^{1}_{s}(f\chi_{E^{c}_{(\cdot)}})\|_{p(\cdot),\gamma_{d}}\leq C_{d,p}$.

Now, if $b>0$ by (\ref{acotacionMb+}) and for all $f\in L^{p(\cdot)}(\mathbb{R}^{d},\gamma_{d})$ with $\|f\|_{p,\gamma_{d}}=1$, we have that
\begin{eqnarray*}
  II &=& \int_{\mathbb{R}^{d}}\left(\int_{B^{c}(x)\cap E_{x}}M(s,x,y)f(y)dy\right)^{p(x)}\gamma_{d}(dx) \\
   &\lesssim&
   \frac{1}{\pi^{d/2}}  \int_{\mathbb{R}^{d}}\left(\int_{B^{c}(x)\cap E_{x}}\frac{e^{-u_{0}}e^{|y|^{2}/p(y)}e^{-|x|^{2}/p(x)}}{t_{0}^{d/2}}f(y)e^{-|y|^{2}/p(y)}dy\right)^{p(x)}dx.
\end{eqnarray*}
Since $p(\cdot)\in\mathcal{P}_{\gamma_{d}}^{\infty}(\mathbb{R}^{d})$, we obtain that $e^{|y|^{2}/p(y)-|x|^{2}/p(x)}\approx e^{(|y|^{2}-|x|^{2})/p_{\infty}}.$ Now by the Cauchy-Schwartz inequality we have,
 $$\left| |y|^{2}-|x|^{2}\right|\leq|x+y| |x-y|,$$
 for all $x,y\in\mathbb{R}^{d}$.\\
  On the other hand, for $b>0$, $|x+y| |x-y|\geq d$ wherever $y\in B_{h}^{c}(x)$.\\
In fact, since $b>0$
\begin{eqnarray}
|x+y|&\geq&|x-y|\label{desnorm1}\\
|x+y|&>&|x|.\label{desnorm2}
\end{eqnarray}
Now, since $y\in B_{h}^{c}(x)$
\begin{itemize}
\item[\underline{Case 1:}] If $|x|\leq1,$ applying (\ref{desnorm1}), we obtain that
$
|x-y|\geq d\left(1\wedge\frac{1}{|x|}\right)=d$ and then 
$$|x-y||x+y|\geq d^{2}\geq d.$$

\item[\underline{Case 2:}] If $|x|>1,$ applying (\ref{desnorm2}), we obtain that
$
|x-y|\geq d\left(1\wedge\frac{1}{|x|}\right)=\frac{d}{|x|}$
and then $$|x-y||x+y|\geq |x-y||x|\geq d.$$
\end{itemize}
Moreover, $t_{0}\approx|x+y| |x-y|/(|x|^{2}+|y|^{2})$. Since $|x|^{2}+|y|^{2}=a<a+b=|x+y|^{2}$, we have that
$$
  t_{0}\geq c \frac{|x+y| |x-y|}{|x|^{2}+|y|^{2}}\geq c\frac{d}{|x+y|^{2}}$$
thus
$$   \frac{1}{t_{0}^{d/2}}\lesssim|x+y|^{d}.$$
Therefore,
\begin{eqnarray*}
  &&\int_{B^{c}(x)\cap E_{x}}\frac{e^{-u_{0}}e^{|y|^{2}/p(y)}e^{-|x|^{2}/p(x)}}{t_{0}^{d/2}}f(y)e^{-|y|^{2}/p(y)}dy\\
   && \hspace{3cm} \lesssim \int_{B^{c}(x)\cap E_{x}}P(x,y)f(y)e^{-|y|^{2}/p(y)}dy
\end{eqnarray*}
where
\begin{eqnarray*}
  P(x,y)=|x+y|^{d}e^{-\alpha_{\infty}|x+y||x-y|} &\hbox{and}&  \alpha_{\infty}=\left(\frac{1}{2}-\left|\frac{1}{p_{\infty}}-\frac{1}{2}\right|\right)>0.
\end{eqnarray*}
It can be proved that $P(x,y)$ is integrable on each variable (for details see \cite{TesSon}) and the value of each integral is independent on $x$ and $y$.\\

 Set $A_{x}=\left\{y:\;\frac{d}{|x|}<|y-x|<\frac{1}{2}\right\}$ and $C_{x}=B^{c}(x,1/2)=\left\{y:\;|y-x|>\frac{1}{2}\right\}$, in consequence $B^{c}(x)\subset A_{x}\cup C_{x}$. Define
\begin{eqnarray*}
  J_{1}=\int_{A_{x}\cap E_{x}}P(x,y)f(y)e^{-|y|^{2}/p(y)}dy &\;\;\hbox{and}\;\;& J_{2}=\int_{C_{x}\cap E_{x}}P(x,y)f(y)e^{-|y|^{2}/p(y)}dy.
\end{eqnarray*}
We will estimate $J_{1}$ first. Observe that, if $y\in A_{x}$, 
 $\frac{3}{4}|x|\leq|y|\leq\frac{5}{4}|x|$ and then $|x|\approx|y|$ hence $|x|\approx|x+y|$, and thus
  \begin{eqnarray*}
    J_{1}
    &\lesssim&\int_{\frac{d}{|x|}<|x-y|}|x|^{d}e^{-\alpha_{\infty}|x||x-y|}f(y)e^{-|y|^{2}/p(y)}dy\\
     &\lesssim& M_{H-L}(fe^{-|\cdot|^{2}/p(\cdot)})(x).
  \end{eqnarray*}
From the hypothesis on $p(\cdot)$ we get
$$
\|M_{H-L}(fe^{-|\cdot|^{2}/p(\cdot)})\|_{p(\cdot)} \lesssim\|fe^{-|\cdot|^{2}/p(\cdot)}\|_{p(\cdot)}\approx\|f\|_{p(\cdot),\gamma_{d}}=1, $$
 and then 
\begin{equation}\label{modmaximal}
\rho_{p(\cdot)}\left(M_{H-L}(fe^{-|\cdot|^{2}
  /p(\cdot)})\right)\lesssim 1.
\end{equation}

To estimate $J_{2}$, we have
 $$ J_{2}\leq\|P(x,\cdot)\chi_{C_{x}}\|_{p^{'}(\cdot)}\leq C,$$
 for details see \cite{DalSco}.This implies that there exists a constant independent on $x$ such that,
$$
  J_{2} = \int_{C_{x}\cap E_{x}}P(x,y)f(y)e^{-|y|^{2}/p(y)}dy\leq C,$$
  thus $$\frac{1}{C}\int_{C_{x}\cap E_{x}}P(x,y)f(y)e^{-|y|^{2}/p(y)}dy \leq 1.$$
We set $g(y)=f(y)e^{-|y|^{2}/p(y)}=g_{1}(y)+g_{2}(y)$, where $g_{1}=g\chi_{\{g\geq1\}}$ y $g_{2}=g\chi_{\{g<1\}}$, applying (\ref{modmaximal}), we have
\begin{eqnarray*}
  II&\lesssim& \int_{\mathbb{R}^{d}}\left(\int_{B^{c}(x)\cap E_{x}}P(x,y)f(y)e^{-|y|^{2}/p(y)}dy\right)^{p(x)}dx \\
   &\lesssim&\int_{\mathbb{R}^{d}}\left(J_{1}\right)^{p(x)}dx+\int_{\mathbb{R}^{d}}\left(J_{2}\right)^{p(x)}dx \\
   &\lesssim&\rho_{p(\cdot)}\left(M_{H-L}(fe^{-|\cdot|^{2}/p(\cdot)})\right) + \int_{\mathbb{R}^{d}}\left(\frac{1}{C}\int_{C_{x}\cap E_{x}}P(x,y)g_{1}(y)dy\right)^{p(x)}dx\\
   &&+\int_{\mathbb{R}^{d}}\left(\frac{1}{C}\int_{C_{x}\cap E_{x}}P(x,y)g_{2}(y)dy\right)^{p(x)}dx \\
  &\lesssim&1+II_{1}+II_{2}
\end{eqnarray*}
Now, we study the integrals $II_{1}$ y $II_{2}$.
\begin{eqnarray*}
  II_{1} &=&  \int_{\mathbb{R}^{d}}\left(\frac{1}{C}\int_{C_{x}\cap E_{x}}P(x,y)g_{1}(y)dy\right)^{p(x)}dx\leq\int_{\mathbb{R}^{d}}\left(\int_{C_{x}\cap E_{x}}P(x,y)g_{1}(y)dy\right)^{p_{-}}dx
\end{eqnarray*}
On the other hand, using Lemma \ref{lema3.26CU} with $G(x)=\frac{1}{C}\int_{C_{x}\cap E_{x}}P(x,y)g_{2}(y)dy\leq1$ and applying the inequality \ref{3.26.1}, we obtain that
\begin{eqnarray*}
  II_{2}&=& \int_{\mathbb{R}^{d}}\left(\int_{C_{x}\cap E_{x}}\frac{1}{C}P(x,y)g_{2}(y)dy\right)^{p(x)}dx =\int_{\mathbb{R}^{d}}(G(x))^{p(x)}dx \\
  &\lesssim& \int_{\mathbb{R}^{d}}(G(x))^{p_{\infty}}dx + \int_{\mathbb{R}^{d}}\frac{dx}{(e+|x|)^{-dp_{-}}}\\
  &=&\int_{\mathbb{R}^{d}}\left(\int_{C_{x}\cap E_{x}}P(x,y)g_{2}(y)dy\right)^{p_{\infty}}+C_{d,p}
\end{eqnarray*}
therefore
\begin{eqnarray*}
  II &\lesssim& \int_{\mathbb{R}^{d}}\left(\int_{C_{x}\cap E_{x}}P(x,y)g_{1}(y)dy\right)^{p_{-}}dx+\int_{\mathbb{R}^{d}}\left(\int_{C_{x}\cap E_{x}}P(x,y)g_{2}(y)dy\right)^{p_{\infty}}dx+C_{d,p}
\end{eqnarray*}
Now, to estimate the last two integrals, we apply H\"{o}lder's inequality.
\begin{eqnarray*}
  &&\int_{\mathbb{R}^{d}}\left(\int_{C_{x}\cap E_{x}}P(x,y)g_{1}(y)dy\right)^{p_{-}}dx\leq\int_{\mathbb{R}^{d}}\left(\int_{\mathbb{R}^{d}}P(x,y)^{\frac{1}{p_{-}^{'}}}P(x,y)^{\frac{1}{p_{-}}}g_{1}(y)dy\right)^{p_{-}}dx\\
&\leq&\int_{\mathbb{R}^{d}}\left(\int_{\mathbb{R}^{d}}(P(x,y))^{p_{-}^{'}/p_{-}^{'}}dy\right)^{p_{-}/p_{-}^{'}}\left(\int_{\mathbb{R}^{d}}(P(x,y))^{p_{-}/p_{-}}g_{1}^{p_{-}}(y)dy\right)^{p_{-}/p_{-}}dx\\
&=&\int_{\mathbb{R}^{d}}\left(\int_{\mathbb{R}^{d}}P(x,y)dy\right)^{p_{-}/p_{-}^{'}}\left(\int_{\mathbb{R}^{d}}P(x,y)g_{1}^{p_{-}}(y)dy\right)dx\\
&\lesssim&\int_{\mathbb{R}^{d}}\int_{\mathbb{R}^{d}}P(x,y)g_{1}^{p_{-}}(y)dydx
\end{eqnarray*}
then, by Fubbini's theorem we get,
\begin{eqnarray*}
  \int_{\mathbb{R}^{d}}\left(\int_{C_{x}\cap E_{x}}P(x,y)g_{1}(y)dy\right)^{p_{-}}dx &\lesssim&\int_{\mathbb{R}^{d}}\int_{\mathbb{R}^{d}}P(x,y)g_{1}^{p_{-}}(y)dydx\\
  &=&\int_{\mathbb{R}^{d}}g_{1}^{p_{-}}(y)\left(\int_{\mathbb{R}^{d}}P(x,y)dx\right)dy\\
  &\lesssim& \int_{\mathbb{R}^{d}}\left(g_{1}(y)\right)^{p(y)}dy \\
  &\lesssim& \int_{\mathbb{R}^{d}}f(y)^{p(y)}e^{-|y|^{2}}dy\lesssim\rho_{p(\cdot),\gamma_{d}}(f)
\end{eqnarray*}
To estimate the integral $\int_{\mathbb{R}^{d}}\left(\int_{C_{x}\cap E_{x}}P(x,y)g_{2}(y)dy\right)^{p_{\infty}}dx$, we proceed in analogous way, but applying the H\"{o}lder's inequality to the exponent $p_{\infty}$, and applying the inequality (\ref{3.26.2}) in Lemma \ref{lema3.26CU}. In consequence we obtain
\begin{eqnarray*}
  \int_{\mathbb{R}^{d}}\left(\int_{C_{x}\cap E_{x}}P(x,y)g_{2}(y)dy\right)^{p_{\infty}}dx &\lesssim&\int_{\mathbb{R}^{d}}g_{2}^{p(y)}(y)dy +C_{d,p} \\
&\lesssim&\rho_{p(\cdot),\gamma_{d}}(f)+C
\end{eqnarray*}
therefore,
\begin{eqnarray*}
  II &\lesssim& \int_{\mathbb{R}^{d}}\left(\int_{C_{x}\cap E_{x}}P(x,y)g_{1}(y)dy\right)^{p_{-}}dx+\int_{\mathbb{R}^{d}}\left(\int_{C_{x}\cap E_{x}}P(x,y)g_{2}(y)dy\right)^{p_{\infty}}dx+C_{d,p}\\
  &\leq&2\rho_{p(\cdot),\gamma_{d}}(f)+C_{d,p}
\end{eqnarray*}
With this we obtain that $\|T^{1}_{s}(f\chi_{E_{(\cdot)}})\|_{p(\cdot),\gamma_{d}}\leq C_{d,p}$, then by homogenity of the norm the result holds for all function $f\in L^{p(\cdot)}(\mathbb{R}^{d},\gamma_{d})$.\\

Hence
\begin{eqnarray*}
  \|T^{1}_{s}f\|_{p(\cdot),\gamma_{d}} &\lesssim& \|T^{1}_{s}(f\chi_{E_{(\cdot)}})\|_{p(\cdot),\gamma_{d}}+\|T^{1}_{s}(f\chi_{E^{c}_{(\cdot)}})\|_{p(\cdot),\gamma_{d}} \lesssim \|f\|_{p(\cdot),\gamma_{d}}. \qed
\end{eqnarray*}
\end{dem}

The proof of boundedness of Ornstein-Uhlenbeck semigroup, Theorem \ref{acotacionsemigrupoOU}, is a immediate consequence of Theorems \ref{acotacionsemigrupoOUlocal} and \ref{acotacionsemigrupoOUglobal}, since for $t>0,$ we have 
\begin{eqnarray*}
  ||T_{t}f||_{p(\cdot),\gamma_{d}}&\leq&\|T^{0}_{s}f\|_{p(\cdot),\gamma_{d}}+\|T^{1}_{s}f\|_{p(\cdot),\gamma_{d}}\leq C\|f\|_{p(\cdot),\gamma_{d}}. \qed\\
\end{eqnarray*}
Additionally, we have,
\begin{teo}\label{continuidadfuerte}
Let $p(\cdot)\in\mathcal{P}_{\gamma_{d}}^{\infty}(\mathbb{R}^{d})\cap LH_{0}(\mathbb{R}^{d})$ with $1<p_{-}\leq p_{+}<\infty,$ and $f\in L^{p(\cdot)}(\mathbb{R}^{d},\gamma_{d})$. The application $t\rightarrow T_{t}f$ is continuous from $[0,\infty)$ to $L^{p(\cdot)}(\mathbb{R}^{d},\gamma_{d})$.
\end{teo}

\begin{dem}
We have to prove that $T_{t}f\rightarrow T_{t_{0}}f$ on $L^{p(\cdot)}(\mathbb{R}^{d},\gamma_{d})$ if $t\rightarrow t_{0}$. By the property of semigroup, it is enough to prove that $T_{t}f\rightarrow f$ in $L^{p(\cdot)}(\mathbb{R}^{d},\gamma_{d})$ if $t\rightarrow0^{+}$.\\

As $f\in L^{p(\cdot)}(\mathbb{R}^{d},\gamma_{d})$, then $f(x)<\infty\;\;a.e.\;x\in\mathbb{R}^{d}$ and $f\in L^{1}(\mathbb{R}^{d},\gamma_{d})$. Let $f_{t}(x)=\left|T_{t}f(x)-f(x)\right|^{p(x)}$, from the pointwise convergence of the Ornstein-Uhlenbeck semigroup (see \cite{TesEbner}), we have,
\begin{eqnarray*}
\lim_{t\rightarrow0^{+}}f_{t}(x)&=&\lim_{t\rightarrow0^{+}}\left|T_{t}f(x)-f(x)\right|^{p(x)}=0,\;a.e.\;\;x\in\mathbb{R}^{d}
\end{eqnarray*}
On the other hand,

\begin{eqnarray*}
\left|T_{t}f(x)-f(x)\right|^{p(x)}\leq2^{p_{+}}\left(\left|T_{t}f(x)\right|^{p(x)}+|f(x)|^{p(x)}\right)&\leq&2^{p_{+}}\left(\left|T^{*}f(x)\right|^{p(x)}+|f(x)|^{p(x)}\right)
\end{eqnarray*}

set $g(x)=2^{p_{+}}\left(\left|T^{*}f(x)\right|^{p(x)}+|f(x)|^{p(x)}\right)\;\;\forall\;x\in\mathbb{R}^{d}$. Then $g$ is integrable, in fact

\begin{eqnarray*}
\int_{\mathbb{R}^{d}}g(x)\gamma_{d}(dx)&=&\int_{\mathbb{R}^{d}}2^{p_{+}}\left(\left|T^{*}f(x)\right|^{p(x)}+|f(x)|^{p(x)}\right)\gamma_{d}(dx)\\
&=&2^{p_{+}}\left(\int_{\mathbb{R}^{d}}\left|T^{*}f(x)\right|^{p(x)}\gamma_{d}(dx)+\int_{\mathbb{R}^{d}}|f(x)|^{p(x)}\gamma_{d}(dx)\right)\\
&=&2^{p_{+}}\left(\rho_{p(\cdot),\gamma_{d}}(T^{*}f)+\rho_{p(\cdot),\gamma_{d}}(f)\right)<\infty
\end{eqnarray*}
since $f$ and $T^{*}f\in L^{p(\cdot)}(\mathbb{R}^{d},\gamma_{d})$.

Applying Lebesgue dominated convergence theorem, we have
\begin{eqnarray*}
\lim_{t\rightarrow0^{+}}\int_{\mathbb{R}^{d}}f_{t}(x)\gamma_{d}(dx)&=&\int_{\mathbb{R}^{d}}\lim_{t\rightarrow0^{+}}f_{t}(x)\gamma_{d}(dx)=0.
\end{eqnarray*}
Thus,
$$
0=\lim_{t\rightarrow0^{+}}\int_{\mathbb{R}^{d}}f_{t}(x)\gamma_{d}(dx)
=\lim_{t\rightarrow0^{+}}\int_{\mathbb{R}^{d}}\left|T_{t}f-f\right|^{p(x)}\gamma_{d}(dx)
=\lim_{t\rightarrow0^{+}}\rho_{p(\cdot),\gamma_{d}}(T_{t}f-f)$$
Then, $\rho_{p(\cdot),\gamma_{d}}(T_{t}f-f)\rightarrow0, \;\;t\rightarrow0^{+}$ and hence $\|T_{t}f-f\|_{p(\cdot),\gamma_{d}}\rightarrow0, \;\;t\rightarrow0^{+}.$
Therefore, $T_{t}f\rightarrow f$ on $L^{p(\cdot)}(\mathbb{R}^{d},\gamma_{d})$ as $t\rightarrow0$.\qed\\
\end{dem}

\subsection{Consequences of the Boundedness of the Ornstein-Uhlenbeck semigroup}

The first consequence of Theorem \ref{acotacionsemigrupoOU} is the proof the boundedness of Poisson-Hermite semigroup, Theorem \ref{acotacionsemigrupoPH}.

\begin{dem}
Let $f\in L^{p(\cdot)}(\mathbb{R}^{d},\gamma_{d})$ with $\|f\|_{p(\cdot),\gamma_{d}}\leq1$, then by Theorem \ref{acotacionsemigrupoOU}, we have for every $s>0$
$$
  ||T_{s}f||_{p(\cdot),\gamma_{d}}\leq C||f||_{p(\cdot),\gamma_{d}}\leq C,$$
    and therefore
 $$ \left\|\frac{T_{s}f}{C}\right\|_{p(\cdot),\gamma_{d}}\leq1.$$
 Thus
 $$\rho_{_{p(\cdot),\gamma_{d}}}\left(\frac{T_{s}f}{C}\right)\leq1.$$

For fixed $t>0$, since the measure $\mu_{t}^{1/2}(ds)$ is a probability measure, using the Jensen's inequality, and Fubini's theorem we get that the modular is less or equal to 1. In fact,
\begin{eqnarray*}
  \rho_{_{p(\cdot),\gamma_{d}}}\left(\frac{P_{t}f}{C}\right) &=& \int_{\mathbb{R}^{d}}\left(\frac{P_{t}f(x)}{C}\right)^{p(x)}\gamma_{d}(dx) 
   \leq \int_{\mathbb{R}^{d}}\int_{0}^{+\infty}\left|\frac{T_{s}f(x)}{C}\right|^{p(x)}\mu_{t}^{1/2}(ds)\gamma_{d}(dx) \\
   &=&\int_{0}^{+\infty}\int_{\mathbb{R}^{d}}\left|\frac{T_{s}f(x)}{C}\right|^{p(x)}\gamma_{d}(dx)\mu_{t}^{1/2}(ds) \\
   &=& \int_{0}^{+\infty}\rho_{_{p(\cdot),\gamma_{d}}}\left(\frac{T_{s}f}{C}\right)\mu_{t}^{1/2}(ds)\leq1.
\end{eqnarray*}

Thus, $P_{t}f\in L^{p(\cdot)}(\mathbb{R}^{d},\gamma_{d})$ and $\|P_{t}f\|_{p(\cdot),\gamma_{d}}\leq C,\;\;\;\forall\;\;t>0$.\\

Now, by homogeneity of the norm and the linearity of $P_{t}$ we obtain the general result.
$$  \|P_{t}f\|_{p(\cdot),\gamma_{d}} \leq C \|f\|_{p(\cdot),\gamma_{d}},$$
for any function $f\in L^{p(\cdot)}(\mathbb{R}^{d},\gamma_{d})$ and $t>0$.\qed \\
\end{dem}

 Additionally, as consequence of Theorem \ref{acotacionsemigrupoPH} we obtain the boundedness of Gaussian Bessel potentials, Theorem \ref{acotacionpotencialBessel}.

\begin{dem}
Let $f\in L^{p(\cdot)}(\mathbb{R}^{d},\gamma_{d})$ with $\|f\|_{p(\cdot),\gamma_{d}}\leq1$, we already know, from the proof of  Theorem \ref{acotacionsemigrupoPH} that, for every $s>0$,
$
  ||P_{s}f||_{p(\cdot),\gamma_{d}}\leq C||f||_{p(\cdot),\gamma_{d}}\leq C
  $
  and therefore
  $\rho_{_{p(\cdot),\gamma_{d}}}\left(\frac{P_{s}f}{C}\right)\leq1.$

Now, for fixed $\beta>0$, using the Jensen's inequality and Fubini's theorem, we get,
\begin{eqnarray*}
  \rho_{_{p(\cdot),\gamma_{d}}}\left(\frac{\mathscr{J}_{\beta}f}{C}\right) &=& \int_{\mathbb{R}^{d}}\left|\frac{\mathscr{J}_{\beta}f(x)}{C}\right|^{p(x)}\gamma_{d}(dx) \\
   &\leq& \int_{\mathbb{R}^{d}}\frac{1}{\Gamma(\beta)}\int_{0}^{+\infty}s^{\beta-1}e^{-s}\left|\frac{P_{s}f(x)}{C}\right|^{p(x)}ds\;\gamma_{d}(dx) \\
   &=& \frac{1}{\Gamma(\beta)}\int_{0}^{+\infty}s^{\beta-1}e^{-s}\int_{\mathbb{R}^{d}}\left|\frac{P_{s}f(x)}{C}\right|^{p(x)}\gamma_{d}(dx)\;ds \\ \\
   &=& \frac{1}{\Gamma(\beta)}\int_{0}^{+\infty}s^{\beta-1}e^{-s}\rho_{_{p(\cdot),\gamma_{d}}}\left(\frac{P_{s}f}{C}\right)ds\leq1.
\end{eqnarray*}

Thus$\mathscr{J}_{\beta}f\in L^{p(\cdot)}(\mathbb{R}^{d},\gamma_{d})$ and 
$$\|\mathscr{J}_{\beta}f\|_{p(\cdot),\gamma_{d}}\leq C,$$
for any $\beta>0$.
Now, again by homogeneity of the norm and linearity of $\mathscr{J}_{\beta}$ we get the general result, 
\begin{eqnarray*}
  \|\mathscr{J}_{\beta}f\|_{p(\cdot),\gamma_{d}} &\leq&C \|f\|_{p(\cdot),\gamma_{d}}
\end{eqnarray*}
for any function $f\in L^{p(\cdot)}(\mathbb{R}^{d},\gamma_{d}).$ \qed
\end{dem}


\end{document}